\documentclass{amsart}
\usepackage{times,amssymb,amsmath,amsfonts,float,nicefrac,subfigure,float}
\usepackage{euscript,graphics}
\usepackage[dvips]{epsfig}

\evensidemargin \oddsidemargin
\newtheorem{theorem}{Theorem}[section]

\newcommand{\remove}[1]{}

\newcommand\nc\newcommand
\nc\bfa{{\boldsymbol a}}\nc\bfA{{\bf A}}\nc\cA{{\mathcal A}}
\nc\bfb{{\boldsymbol b}}\nc\bfB{{\bf B}}\nc\cB{{\mathcal B}}
\nc\bfc{{\boldsymbol c}}\nc\bfC{{\bf C}}\nc\cC{{\mathcal C}}
\nc\bfd{{\boldsymbol d}}\nc\bfD{{\bf D}}\nc\cD{{\mathcal D}}
\nc\bfe{{\boldsymbol e}}\nc\bfE{{\bf E}}\nc\cE{{\mathcal E}}
\nc\bff{{\boldsymbol f}}\nc\bfF{{\bf F}}\nc\cF{{\mathcal F}}
\nc\bfg{{\boldsymbol g}}\nc\bfG{{\bf G}}\nc\cG{{\mathcal G}}
\nc\bfh{{\boldsymbol h}}\nc\bfH{{\bf H}}\nc\cH{{\mathcal H}}
\nc\bfi{{\boldsymbol i}}\nc\bfI{{\bf I}}\nc\cI{{\mathcal I}}
\nc\bfj{{\boldsymbol j}}\nc\bfJ{{\bf J}}\nc\cJ{{\mathcal J}}
\nc\bfk{{\boldsymbol k}}\nc\bfK{{\bf K}}\nc\cK{{\mathcal K}}
\nc\bfl{{\boldsymbol l}}\nc\bfL{{\bf L}}\nc\cL{{\mathcal L}}
\nc\bfm{{\boldsymbol m}}\nc\bfM{{\bf M}}\nc\cM{{\mathcal M}}
\nc\bfn{{\boldsymbol n}}\nc\bfN{{\bf N}}\nc\cN{{\mathcal N}}
\nc\bfo{{\boldsymbol o}}\nc\bfO{{\bf O}}\nc\cO{{\mathcal O}}
\nc\bfp{{\boldsymbol p}}\nc\bfP{{\bf P}}\nc\cP{{\mathcal P}}
\nc\bfq{{\boldsymbol q}}\nc\bfQ{{\bf Q}}\nc\cQ{{\mathcal Q}}
\nc\bfr{{\boldsymbol r}}\nc\bfR{{\bf R}}\nc\cR{{\mathcal R}}
\nc\bfs{{\boldsymbol s}}\nc\bfS{{\bf S}}\nc\cS{{\mathcal S}}
\nc\bft{{\boldsymbol t}}\nc\bfT{{\bf T}}\nc\cT{{\mathcal T}}
\nc\bfu{{\boldsymbol u}}\nc\bfU{{\bf U}}\nc\cU{{\mathcal U}}
\nc\bfv{{\boldsymbol v}}\nc\bfV{{\bf V}}\nc\cV{{\mathcal V}}
\nc\bfw{{\boldsymbol w}}\nc\bfW{{\bf W}}\nc\cW{{\mathcal W}}
\nc\bfx{{\boldsymbol x}}\nc\bfX{{\bf X}}\nc\cX{{\mathcal X}}
\nc\bfy{{\boldsymbol y}}\nc\bfY{{\bf Y}}\nc\cY{{\mathcal Y}}
\nc\bfz{{\boldsymbol z}}\nc\bfZ{{\bf Z}}\nc\cZ{{\mathcal Z}}
\nc\od{{\bar d}}\nc\ow{{\bar w}}\nc\odelta{{\bar\delta}}
\nc\ox{{\bar x}}\nc\oy{{\bar y}}\nc\ou{{\bar u}}
\nc\oh{{\bar h}}

\newcommand\reals{{\mathbb R}}
\newcommand\complexes{{\mathbb C}}

\nc\dgv{\delta_{\text{\rm GV}}}
\nc\dcrit{\delta_{\text{\rm{crit}}}}
\nc\Esp{E_{\text{\rm sp}}}
\renewcommand\epsilon{\varepsilon}

\nc\ellone{{\ell_1}}
\nc\elltwo{{\ell_2}}
\nc\ellinf{{{\ell_\infty}}}
\nc\ip[2]{\langle #1,#2\rangle}

\newcommand{\beeq}{\begin{eqnarray*}}
\newcommand{\eneq}{\end{eqnarray*}}

\newcommand{\half}{\nicefrac12}

\begin{document}
\title[Spherical two-distance sets]
{New bounds for spherical two-distance sets}
\thanks{{\em Date}\/: December 29, 2012.\/ Research
supported in part by NSF grants DMS1101697,
CCF0916919, and CCF0830699 and NSA grant H98230-101013.}%

\author[A. Barg]{Alexander Barg$^\ast$}\thanks{$^\ast$
Dept. of Electrical and Computer Engineering and Institute for Systems
Research, University of Maryland, College Park, MD 20742,
and Institute for Information Transmission Problems,
Russian Academy of Sciences, Moscow, Russia. Email: abarg@umd.edu.}
\author[W.-H. Yu]{Wei-Hsuan Yu$^\dag$} \thanks{$^\dag$
Dept. of Mathematics and Institute for Systems Research, University
of Maryland, College Park, MD 20742. Email: mathyu@math.umd.edu}

\begin{abstract}
A spherical two-distance set is a finite collection of unit
vectors in $\reals^n$ such that the distances between any
two distinct vectors assume only two values. We use the semidefinite
programming method to compute improved estimates of the maximum
size of spherical two-distance sets. Exact answers are found for
dimensions $n=23$ and $40\le n\le 93\; (n\ne 46,78)$ where previous results gave
divergent bounds.
\end{abstract}
\maketitle
\section{Introduction}
This paper is devoted to the application of the semidefinite
programming method to estimates of the size of the largest
possible two-distance set on the sphere $S^{n-1}(\reals)$. A spherical two-distance
set is a finite collection $\cC$ of unit vectors in $\reals^n$
such that the set of distances between any two distinct vectors in
$\cC$ has cardinality two. Estimating the maximum size $g(n)$ of
such a set is a classical problem in distance geometry that has
been studied for several decades.

We begin with an overview of known results. A lower
bound on $g(n)$ is obtained as follows. Let
$e_1,\dots,e_{n+1}$ be the standard basis in $\reals^{n+1}$.
The points $e_i+e_j, i\ne j$ form a spherical two-distance set in the
plane $x_1+\dots+x_{n+1}=2$ (after scaling), and therefore
   \begin{equation}\label{eq:n}
     g(n)\ge n(n+1)/2, \quad n\ge 2.
  \end{equation}

The first major result for upper bounds was obtained by
Delsarte, Goethals, and Seidel \cite{del77b}. They proved that, irrespective of the actual
values of the distances, the following ``harmonic'' bound holds true:
  \begin{equation}\label{eq:dgs}
   g(n)\le n(n+3)/2.
  \end{equation}
They also showed that this bound is tight for dimensions $n=2,6,22$ where it
is attained by sets of equiangular lines. Moreover, the results of \cite{del77b},
Bannai et al.~\cite{ban04}, and Nebe and Venkov \cite{neb12}
imply that $g(n)$ can attain the harmonic bound
only if $n=(2m+1)^2-3, m\ge 1$ with the exception of an infinite sequence of values of
$m$ that begins with $m=3,4,6,10,12,22,38,30,34,42,46$. Therefore, unless $n$
is of the above form, $g(n)\le n(n+3)/2-1.$
These results are proved using the link between 2-distance sets and tight spherical
4-designs established in \cite{del77b}. 

Another advance in estimating the function $g(n)$ was made by Musin
\cite{mus09a}. Let $\cC=\{z_1,z_2,\dots\}$ and suppose that
$z_i\cdot z_j\in\{a,b\}, i\ne j,$
where $2-2a, 2-2b$ are the values of the squared distances between the points.
Musin proved that
   \begin{equation}\label{eq:m}
    |\cC|\le n(n+1)/2 \quad\text{if } a+b \ge 0.
   \end{equation}
He then used Delsarte's linear programming method 
to prove that $g(n)=n(n+1)/2$ if
$7\le n\le 39, n\ne 22,23.$


Here we make another step for spherical two-distance sets, extending the range of
dimensions in which the bound \eqref{eq:n} is tight.  The state of the art for $g(n)$ can be summarized as follows.
\begin{theorem} 
We have $g(2)=5, g(3)=6, g(4)=10, g(5)=16, g(6)=27,g(22)=275,$
   \begin{align}
n(n+1)/2\le &g(n)\le n(n+3)/2-1, \quad n=46,78\label{eq:results2}\\
     &g(n)=n(n+1)/2, \quad 7\le n\le 93,
n\ne 22,46,78,\label{eq:results1}     
   \end{align}
and $4465\le g(94)\le 4492.$ If $n\ge 95,$ then
  $   g(n)\le n(n+3)/2 \text{ or } n(n+3)/2-1$
as detailed in the remarks after Eq. \eqref{eq:dgs} above.
\end{theorem}
The part of this theorem that is established in the present paper relates to dimensions
$n=23$ and $40\le n\le 94, n\ne 46,78.$ Our results are 
computational in nature and are obtained using the semidefinite programming method.
The other parts of this theorem follow from the results in \cite{del77b,ban04,mus09a,neb12}.

As far as actual constructions of spherical two-distance sets are concerned, rather little
is known beyond the set of midpoints of the edges of a regular simplex mentioned above.
Another way of constructing such sets is to start with a set of equiangular lines
in $\reals^n$ \cite{lem73}. If the angle between each pair of lines is $\alpha,$ then
taking one point from each pair of points on $S^{n-1}$ defined by the line, we obtain
a two-distance set with $a=\alpha,b=-\alpha.$ The largest possible number of equiangular
lines in $\reals^n$ is $n(n+1)/2$ (this result is due to Gerzon, see \cite{lem73}). 
This bound is attained for $n=3,7,23.$ For instance, for $n=3$ the
set of 6 lines is obtained from 6 diagonals of the icosahedron, which gives many
ways of constructing inequivalent spherical two-distance sets of cardinality 6. 
The only three instances in which the known spherical two-distance sets are of cardinality
greater than $n(n+1)/2$ occur in dimensions $n=2,6$ and $22.$

\section{Positive definite matrices and SDP bounds}
A semidefinite program is an optimization problem of the form
 \begin{equation}\label{eq:SDP}
   \max\{\langle X,C\rangle|\; X\succeq 0,\; \langle X,A_i\rangle=b_i, i=1,\dots,m\},
 \end{equation}
 where $X$ is an $n\times n$ variable matrix, $A_1,\dots,A_m$ and $C$ are given Hermitian
 matrices, $(b_1,\dots, b_m)$ is a given vector and $\langle X, Y\rangle=\text{trace}\, (Y^\ast X)$ is the
 inner product of two matrices. Semidefinite programming is an extension of linear programming
 that has found a range of applications in combinatorial optimization, control theory, 
 distance geometry, and coding theory. General introduction to semidefinite programming is given,
for instance, 
 in \cite{ben01}. 

The main problem addressed by the SDP method in distance geometry is related to 
deriving bounds on the cardinality of point
sets in a metric space $\cX$ with a given set of properties such as a given minimum separation between distinct points in the set.
The SDP method has its roots in harmonic analysis of the isometry
group of the metric space in question. It is broadly applicable in both finite and
compact infinite spaces. Examples of the former include the Hamming and Johnson spaces,
their $q$-analogs, other metric spaces on the set of $n$-strings over a finite alphabet,
as well as the finite projective space. The main example in the infinite case is given by
real and complex spheres, although the SDP method is also applicable in other compact
homogeneous spaces.
Working out the details in each example is a nontrivial task that includes
analysis of irreducible modules in the space of functions $f: \cX\to \complexes$ 
under the action of the isometry group $G$ of $\cX$.
The zonal matrices that arise in this analysis initially have large size that can be reduced
relying on symmetries arising from the group action. This gives rise to an SDP optimization problem
that is solved by computer for a given set of dimensions (the numerical part is also not straightforward
and rather time-consuming). Foundations and analysis of particular cases have been the subject of
a considerable number of research and overview publications in the last decade; see in particular
recent surveys \cite{bac12,bac09} and references therein.

The origins of the SDP method and the discussed applications can be traced back to the work
of Delsarte \cite{del73a} which introduced the machinery of association schemes in the analysis
of point configurations (codes) in finite spaces. Delsarte derived linear programming (LP) bounds on 
the cardinality of a set of points in the space under the condition on the minimum separation of
distinct points in the set. Delsarte's results were linked to harmonic analysis and group
representations in the works of Delsarte, Goethals and Seidel \cite{del77b} (for the case
$S^{n-1}$) and Kabatyansky and Levenshtein \cite{kab78} (for general compact symmetric spaces).

From now on we focus on the case $\cX=S^{n-1}.$ Let $G_k^{(n)}(t), k=0,1,\dots$ denote the Gegenbauer polynomials of degree $k$.
They are defined recursively as follows: $G_0^{(n)}\equiv 1, G_1^{(n)}(t)=t,$
 and
  $$
    G_k^{(n)}(t)=\frac{(2k+n-4)tG_{k-1}^{(n)}(t)-(k-1)G_{k-2}^{(n)}(t)}
   {k+n-3}, \quad k\ge 2.
  $$
Delsarte et al. \cite{del77b} showed that for any finite set of points $\cC\subset S^{n-1}$
  \begin{equation}\label{eq:Gp}
    \sum_{(x,y)\in \cC^2} G_k^{(n)}(x\cdot y) \ge 0, \quad k=1,2,\dots .
  \end{equation}
The proof of this inequality in \cite{del77b} used the addition formula for spherical harmonics. 
An earlier, geometric
proof of \eqref{eq:Gp} had been given by Schoenberg \cite{sch42}, although his work 
was not known to researchers in the area discussed until at least the 1990s.

Positivity conditions \eqref{eq:Gp} give rise to the LP bound on the cardinality of spherical two-distance sets.
\begin{theorem} \label{thm:DGS}
{\rm (Delsarte et al. \cite{del77b})} Let $\cC\subset S^{n-1}$ be a finite set and suppose that 
$x\cdot y\in\{a,b\}$ for any $x,y\in \cC.$
Then   $$
  |\cC|\le\max\Big\{1+\alpha_1+\alpha_2: \;
             1+\alpha_1G_i^{(n)}(a)+\alpha_2 G_i^{(n)}(b)
\ge 0, i=0,1,\dots,p;\alpha_j\ge 0, j=1,2\Big\}.
  $$
\end{theorem}
In this theorem $\alpha_1,\alpha_2$ are the optimization variables that refer 
to the number of ordered pairs of points in $\cC$ with inner product $a$ and $b$, respectively.
For instance, $
    \alpha_1=|\cC|^{-1}\sharp\{(z_1,z_2)\in \cC^2: z_1\cdot z_2=a\},
  $
This theorem is a specialization of a more general LP bound on spherical codes of
\cite{del77b,kab78}.

Applications of semidefinite programming in coding theory and distance geometry
gained momentum after the pioneering work of Schrijver \cite{sch05} that derived
SDP bounds on codes in the Hamming and Johnson spaces. Schijver's 
approach was based on the so-called Terwilliger algebra of the association scheme
and formed a far-reaching generalization of the work of Delsarte \cite{del73a}. 
Elements of the groundwork for SDP bounds in the Hamming space were
laid by Dunkl \cite{dun76}, although this connection was also made somewhat later \cite{val09a}.
We refer to \cite{mar09} for a detailed general survey of the approach via association schemes
and further references. 

SDP bounds for the real sphere were derived by Bachoc and Vallentin \cite{bac08a} in the context
of the kissing number problem. The kissing number $k(n)$ is the maximum number of unit spheres
that can touch a unit sphere without overlapping, i.e. the maximum number of points on the sphere
such that the angular separation between any pair of them is at least $\pi/3.$ 
Following \cite{bac08a}, define a $(p-k+1)\times(p-k+1)$ matrix $Y_k^n(u,v,t), k\ge 0$
by setting
  $$
    (Y_k^n(u,v,t))_{ij}=u^iv^j ((1-u^2)(1-v^2))^{k/2}
      G_k^{(n-1)}\Big(\frac{t-uv}{\sqrt{(1-u^2)(1-v^2)}}\Big)
  $$
where $p$ is a positive integer, and a matrix $S_k^n(u,v,t)$ by setting
  \begin{equation}\label{eq:S}
    S_k^n(u,v,t)=\frac16\sum_\sigma Y_k^n(\sigma(u,v,t)),
  \end{equation}
where the sum is over all permutations on 3 elements. Note that
$(S_k^n(1,1,1))_{ij}=0$ for all $i,j$ and all $k\ge 1.$ One of the main results of
\cite{bac08a} is that for any finite set of points $\cC\subset S^{n-1}$
  \begin{equation}\label{eq:Sp}
   \sum_{(x,y,z)\in \cC^3} S_k^n(x\cdot y,x\cdot z,y\cdot z)\succeq 0
  \end{equation}
The matrices $S_k^n$ play the role of the constraints $A_i$ in the general SDP problem \eqref{eq:SDP}. 
Positivity constraints \eqref{eq:S} give rise to a general SDP bound
on the cardinality of point sets obtained in \cite{bac08a}, where it was used to improve
upper bounds on $k(n)$ in small dimensions. In the next section
we state a specialization of this bound for the case of 2-distance sets.

As a final remark, we note that constraints \eqref{eq:Gp} arise from the unrestricted
action of $G$ on $S^{n-1}.$ Constraints \eqref{eq:S} are obtained
by considering only actions that fix an arbitrary given point on the sphere.
Further SDP bounds can be obtained by considering zonal matrices that arise from actions that
fix any given number of points; however even for two points, actual evaluation
of the bounds requires significant computational effort \cite{val10}.

\section{The bounds}

The general SDP bound on spherical codes of \cite{bac08a} specializes to our case as follows.
\begin{theorem}\label{thm:SDP} Let $\cC$ be a spherical two-distance set
with inner products $a$ and $b$.
Let $p$ be a positive integer.
The cardinality $|\cC|$ is bounded above
by the solution of the following semidefinite programming problem:
  \begin{equation}\label{eq:sdp1}
  1+\nicefrac13\max (x_1+ x_2)
  \end{equation}
\centerline{\text{\rm subject to }}
  \begin{equation}
\label{eq:sdp2} \begin{pmatrix}
1 & 0  \\
0 & 0\end{pmatrix} +
 \frac{1}{3}\begin{pmatrix}
0 & 1  \\
1 & 1\end{pmatrix} (x_1 + x_2)
+ \begin{pmatrix}
0 & 0  \\
0 & 1
\end{pmatrix}  (x_3+x_4+x_5+x_6)  \succeq 0
  \end{equation}
  \begin{equation}\label{eq:sdp3}
3+G_i^{(n)}(a)x_1+G_i^{(n)}(b)x_2\ge 0, \quad i=1,2,\dots,p
  \end{equation}
   \begin{align}\label{eq:sdp4}
S^n_i(1,1,1)+S^n_i(a,a,1)x_1+&S^n_i(b,b,1)x_2+S^n_i(a,a,a)x_3\\
  &+S^n_i(a,a,b)x_4+S^n_i(a,b,b)x_5+S^n_i(b,b,b)x_6 \succeq 0, \quad i=0,1,
\dots, p  \nonumber  \end{align}
$$x_j\ge 0, j=1,\dots,6,$$
where $S_i(\cdot,\cdot,\cdot)$ are $(p-i+1)\times(p-i+1)$ matrices
defined in \eqref{eq:S}.
\end{theorem}
In this theorem the variables $x_1,x_2$ refer to the number of ordered pairs
of vectors in $\cC$ with inner product $a$ and $b$ respectively; namely we have
$x_i=3\alpha_i, i=1,2.$
We note that the SDP problem seeks to optimize the same linear form as the LP problem, but
adds more constraints on the configuration. Because of this,
Theorem \ref{thm:SDP} usually gives tighter
bounds than Theorem \ref{thm:DGS}. This fact is evident from the table below and is
also known from the
calculation of kissing numbers in \cite{bac08a}.

\subsection{Calculation of the bound} Several remarks are in order.
First, implementation of SDP for two-distance sets differs from earlier computations
in \cite{bac08a,val10} in that in our case there are no limits
on the minimum separation of the points.
Next, we restrict our calculations
to the case $p\le 5$ as no improvement is observed for larger
values. Finally, by a result Larman et al. \cite{lar77}, if $|\cC|\ge
2n+3$ then the inner products $a,b$ are related by
$b=b_k(a)=(ka-1)/(k-1)$ where $k\in\{2,\dots,
\lfloor(1+\sqrt{2n})/2\rfloor\}$ is an integer. Thus we obtain a
family of SDP bounds parametrized by $a.$ Since $b_k(a)\ge -1,
a+b_k(a)<0,$ we get that $a\in I_k:=[0,\frac 1{2k-1}).$ Moreover,
if $-1\le b<a\le0,$ then $|\cC|$ cannot be large by the Rankin
bounds \cite{ran55}, and if $a+b\ge0$ then $|\cC|$ is bounded by
\eqref{eq:m}. We conclude as follows.
\begin{theorem} Let {\rm SDP}$(a)$ be the solution of the SDP problem
\eqref{eq:sdp1}-\eqref{eq:sdp4}, where $b=b_k(a)$. Let $\cC$ be a spherical
two-distance set with inner products $a,b$, then
   $$
 |\cC|\le \begin{cases} n(n+1)/2, &a+b\ge 0\\
    {\rm SDP}(a), &a\in I_k\\
    n+1, &-1\le b<a<0.
  \end{cases}
  $$
\end{theorem}
For instance, for $n=23, k=3$ we obtain that $I_k=[0,0.2)$. Partitioning
$I_k$ into a number of small segments, we plot the value $\text{SDP}(a)$
as a function of $a$ evaluated at the nodes of the partition.
The result is shown in Fig.~\ref{fig:subfig1}. A part of the segment
around the maximum appears in Fig.~\ref{fig:subfig2}.
This computation gives an indication of the answer, but in principle
the value $\text{SDP}(a)$ could oscillate between the nodes of the partition.
Ruling this out requires perturbation analysis of the SDP problem
which is not immediate.

\subsubsection{Dual problem}
The dual problem of \eqref{eq:sdp1}-\eqref{eq:sdp4} has the following
form.
  \begin{equation}\label{eq:dual}
   1+\min\Big\{\sum_{i=1}^p \alpha_i+\beta_{11}+\langle F_0,S_0^n(1,1,1)\rangle\Big\}
  \end{equation}
\centerline{subject to}
  $$
  \begin{pmatrix}\beta_{11}&\beta_{22}\\\beta_{12}&\beta_{22}\end{pmatrix}\succeq0
  $$
\begin{equation}
  2\beta_{12}+\beta_{22}+\sum_{i=1}^p (\alpha_iG_i^{(n)}(a)+3\langle F_i,S_i^n(a,a,1)\rangle )
  \le -1\label{eq:dual1}
\end{equation}
\begin{equation}  2\beta_{12}+\beta_{22}+\sum_{i=1}^p (\alpha_iG_i^{(n)}(b)+
3\langle F_i,S_i^n(b,b,1)\rangle ) \le -1\label{eq:dual2}\end{equation}
 \begin{gather} \beta_{22}+\sum_{i=0}^p\langle F_i,S_i^n(y_1,y_2,y_3)\rangle\le 0\label{eq:dual3}\\
 \hspace*{2in} \text{where }(y_1,y_2,y_3)\in\{(a,a,a),(a,a,b),(a,b,b),(b,b,b)\}\nonumber\\
  \alpha_i\ge 0, \;F_i\succeq0, \;i=1,\dots,p.\nonumber
\end{gather}

We need to estimate from above the maximum value of this problem over $a\in I_k=[a_1,a_2].$
Accounting for a continuous value set of the parameter in SDP problems is a
challenging task.
We approach it by employing the sum-of-squares method.
Constraints \eqref{eq:dual1}-\eqref{eq:dual3} impose positivity conditions on some univariate
polynomials of $a$ for $a\in I_k.$ The following sequence of steps
transforms the constraints to semidefinite conditions. Observe that a polynomial $f(a)$ of degree at most $m$
satisfies $f(a)\ge 0$ for $a\in I_k$
if and only if the polynomial of degree at most $2m$
  $$
f^+(a)=(1+a^2)^m f\Big(\frac {a_1+a_2a^2}{1+a^2}\Big) \ge 0
  $$
for all $a\in \reals.$
Next, a polynomial nonnegative on the entire real axis
can be written as a sum of squares, $f(x)=\sum_i r_i^2(x),$ where the
$r_i$ are polynomials. Further, by a result of Nesterov \cite{nes00},
a polynomial $f(x)$ of degree $2m$ is a sum of squares
if and only if there exists a positive semidefinite matrix $Q$ such that $f=XQX^t,$
where $X=(1,x,x^2,\dots, x^m).$ Thus, constraints \eqref{eq:dual1}-\eqref{eq:dual3}
can be transformed to semidefinite conditions.

As a result, we obtain an SDP problem that can be solved by computer. We solved the resulting problem
for $7\le n\le 96$
using the Matlab toolbox SOSTOOLS \cite{S} in the YALMIP environment \cite{YALMIP}.
An advantage in using SOSTOOLS is that it accepts $a$ as an SDP variable, thereby accounting for
all the values of $a$ in the segment. Thus, we obtain the value $\max\text{SDP(a)}, a\in I_k.$
However, this may impose excessive constraints on the value of the SDP problem because
all the conditions for different values of $a$ are involved at the same time. To work around
this accumulation, we use a sub-partitioning of the segment $I_k$ into smaller
segments. For each of them, SOSTOOLS outputs the largest value of the minimum of the SDP problem
over all $a$ in the segment. It turns out that, in many cases, 
the maximum of these solutions is smaller than $\max\text{SDP(a)}, a\in I_k$ computed directly
by the package.
The estimates of the answer computed from the primal problem serve as a guidance of the 
needed step length of the partition. The solution of the sum-of-squares
SDP optimization problem provides a rigorous proof for
the estimates obtained by discretizing the primal problem \eqref{eq:sdp1}-\eqref{eq:sdp4}.
For instance, for $n=23$ we partition $I_3$ into 20 subsegments, finding 276.5 as the
maximum value of the dual SDP problem for $a \in I_3,$ etc.

\remove{We proceed as in \cite{bac08a}, relying a theorem of M.
Putinar about representation of positive polynomials as sums of
squares and the corresponding relaxation of the SDP problem.
Letting $I_k=[a_1,a_2],$ we observe that $a\in I_k$ is and only if
   $$
   r(a)=-(a-\half(a_1+a_2))^2+(1/4)(a_1-a_2)^2\ge 0.
   $$
    A polynomial
$q(x_1,\dots,x_n)$ of degree $2m$ is a sum of squares (i.e.,
$q(x)=\sum_i r_i^2$) if and only if there is a positive
semidefinite matrix $Q$ such that $q=z^t Q z,$ where
$z=(1,x_1,\dots,x_n,x_1x_2,\dots,x_{n-1}x_n,\dots,x_n^m)$ is the
$\binom{n+m}m$-vector of monomials. We can write inequality
\eqref{eq:dual1} in the form
  $$
  -1-2\beta_{12}-\beta_{22}-\sum_{i=1}^p (\alpha_iG_i^{(n)}(a)+3\langle F_i,S_i^n(a,a,1) )
  \rangle=q(a)+p(a)q_1(a),
  $$
where $q,q_1$ are squares of polynomials of $a$. Similar
inequalities can be written to replace
\eqref{eq:dual2}-\eqref{eq:dual3}.}

\subsubsection{Results} The results of the calculation are summarized
in the table below. The part of the table for $7\le n\le 40,$
except for the values of the SDP bound, is from \cite{mus09a}. The
improvement provided by Theorem \ref{thm:SDP} over the LP bound is
quite substantial even for relatively small dimensions. The LP
bound is above $n(n+1)/2$ for $n\ge 40$ and is not included
starting with $n=41.$ The cases $n=46,78$ and $n \geq 94$ are not
resolved by SDP, although for $n=94$ we still obtain an improvement
over the harmonic bound \eqref{eq:dgs}. 
The value of $k$ shown in the
table accounts for the largest value of the SDP problem among the possible choices of $k$.
This guarantees that the value SDP$(a)$ is equal to or smaller than the number
in the table for all the possible values of the inner products $a,b$ in the point set.

\begin{figure}[H]
\subfigure[The value SDP$(a)$ for $n=23$]{
\includegraphics[scale=0.5]{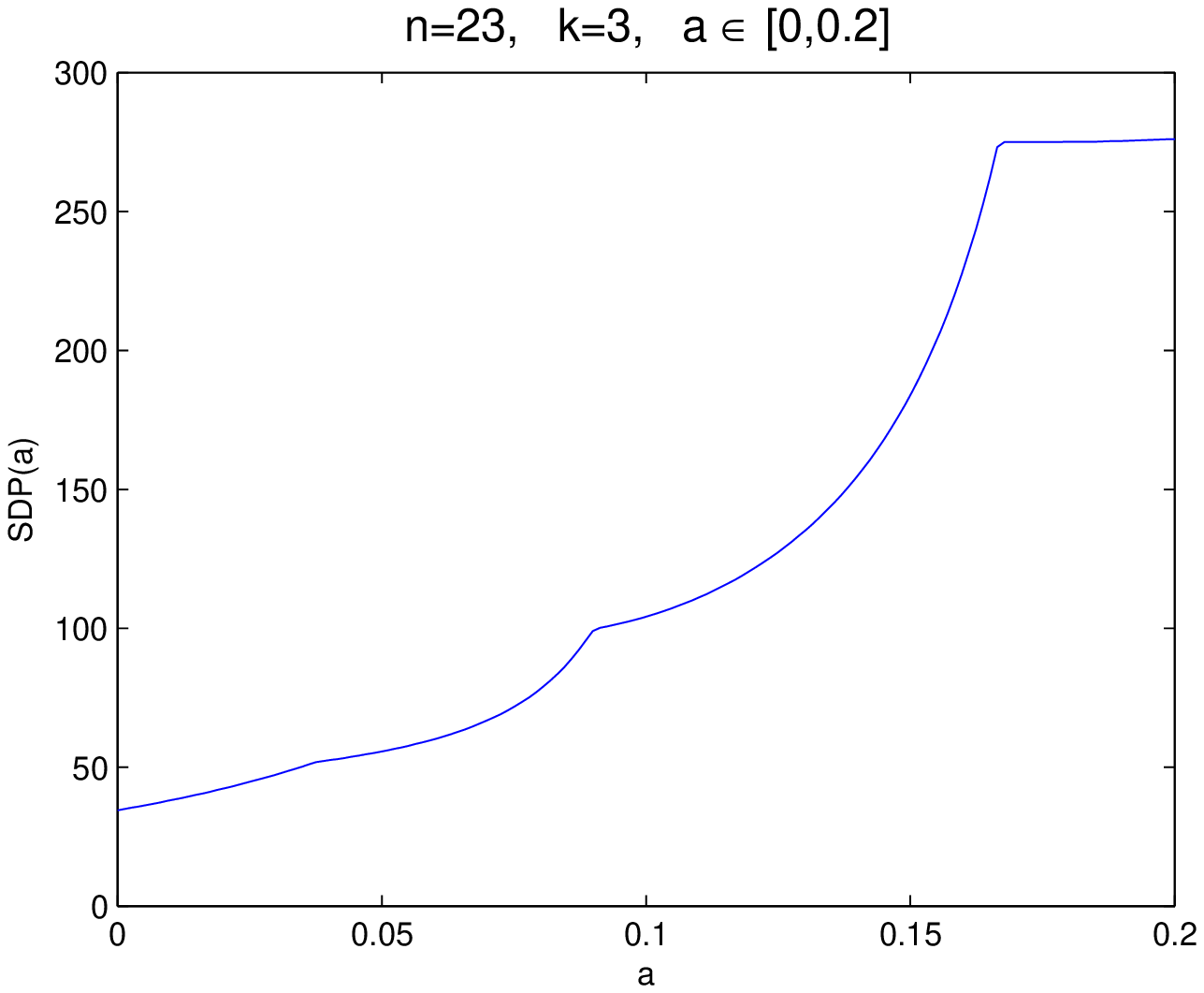}
\label{fig:subfig1}
}
\subfigure[The neighborhood of the maximum]{
\includegraphics[scale=0.5]{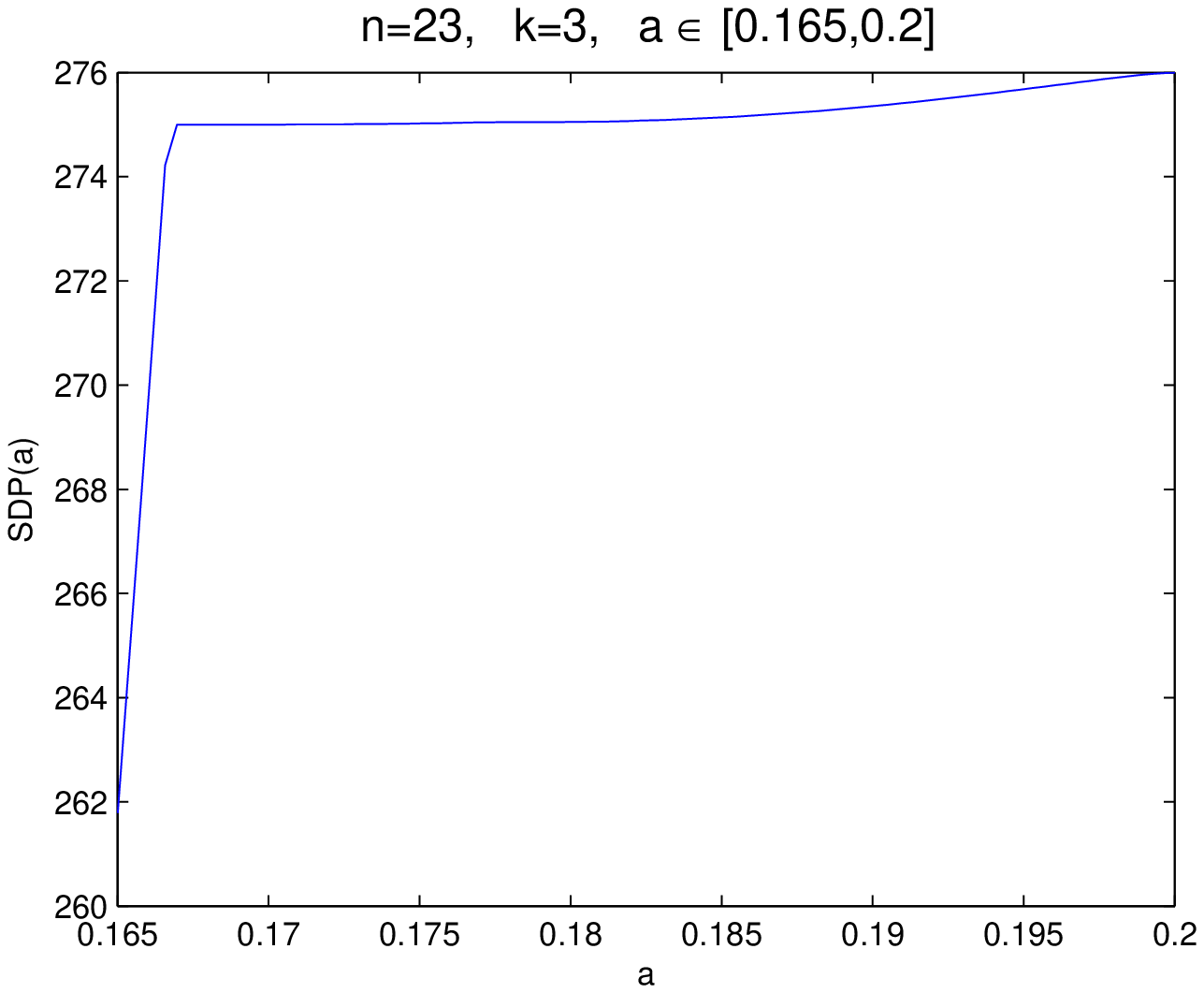}
\label{fig:subfig2}}
\caption{Evaluation of the SDP bound on $g(23)$}\end{figure}

Notice that for $n=46,78$ the SDP bound coincides with the
bound \eqref{eq:dgs}. For $n=23$ the results of
\cite{mus09a} leave two possibilities, $g(n)=276$ and 277. The SDP
bound resolves this for the former, establishing the corresponding part
of the claim in
\eqref{eq:results1}. As is seen from Fig.~\ref{fig:subfig2}, the largest value
of SDP$(a)$ is attained for $a=0.2$ and is equal to 276. This case corresponds to 276
equiangular lines in $R^{23}$ with angle $\arccos 0.2$, which can be constructed
either using strongly regular graphs or the Leech lattice (see \cite{lem73} for details).

\vspace*{.05in}{\em Acknowledgment:} We are grateful to Chao-Wei Chen and
Johan L\"{o}fberg for their significant help with the Matlab implementation.

\remove{In conclusion  note that Musin and Nozaki \cite{mus11} computed bounds on
the size of spherical sets with three and four distances for small dimensions.
It is easy to write an SDP bound on their size similar to Theorem \ref{thm:SDP},
which is likely to improve on the known results.
However, the volume of computations involved in the actual optimization and
sensitivity analysis at this point looks prohibitively large to pursue this
line.}

\providecommand{\bysame}{\leavevmode\hbox to3em{\hrulefill}\thinspace}
\providecommand{\href}[2]{#2}

\begin{table}
\begin{center}{\small Bounds on two-distance sets. The starred rows correspond to dimensions for
which\\ the value of $g(n)$ is not known exactly. }\end{center}

\vspace*{.1in}{\footnotesize\begin{tabular}{|l|c|c|c|c|@{\hspace*{.1in}}|l|c|c|c|}
\hline
$n$ & LP bound & SDP bound & $n(n+1)/2$ & $k$&
                            $n$ &  SDP bound & $n(n+1)/2$ & $k$\\
\hline
7& 28 & 28 & 28 &2         &52& 1128 & 1378  &4\\
8& 31 & 28 & 36 &2         &53& 1128 & 1431  &4\\
9& 34 & 29 & 45 &2         &54& 1128 & 1485  &4\\
10& 37 & 29 & 55 &2        &55& 1128 & 1540  &4\\
11& 40 & 29 & 66 &2        &56& 1128&  1596  &4 \\
12& 44 & 28 & 78 &2        &57& 1162&  1653  &2\\
13& 47 & 29 & 91 &3        &58& 1200 & 1711  &2 \\
14& 52 & 35 & 105 &2       &59& 1240&  1770  &2 \\
15& 56 & 41 &120 &3        &60& 1282&  1830  &2\\
16& 61 & 50 &136 &3        &61& 1324&  1891  &2\\
17& 66 & 60 &153 &3        &62& 1372&  1953  &2\\
18& 76 & 75 &171 &3        &63& 1428 & 2016  &2 \\
19& 96 & 95 &190 &3        &64& 1482&  2080  &2 \\
20& 126 & 124& 210 &3      &65& 1540&  2145  &2 \\
21& 176 & 174 &231 &3      &66& 1604&  2211  &2 \\
22& 275 & 275 &253 &3      &67& 1672&  2278  &2 \\
23& 277 & 276 &276 &3      &68& 1745&  2346  &2 \\
24& 280 & 276 &300 &3      &69& 1822&  2415  &2 \\
25& 284 & 276 &325 &3      &70& 1907&  2485  &2 \\
26& 288 & 276 &351 &3      &71& 1999&  2556  &2 \\
27& 294 & 276 &378 &3      &72& 2097&  2628  &2 \\
28& 299 & 276 &406 &3      &73& 2206&  2701  &2 \\
29& 305 & 276 &435 &3      &74& 2325&  2775  &2 \\
30& 312 & 276 &465 &3      &75& 2394&  2850  &2  \\
31& 319 & 276 &496 &3      &76& 2468&  2926  &2 \\
32& 327 & 276 &528 &3      &77& 2542&  3003  &2 \\
33& 334 & 276 &561 &3      &$78^\ast$& 3159 & 3081  &2 \\
34& 342 & 276 &595 &3      &79& 3160 & 3160  &4 \\
35& 360 & 276 &630 &2      &80& 3160 & 3240  &4\\
36& 416 & 276 &666 &2      &81& 3160 &  3321  &4 \\
37& 488 & 276 &703 &2      &82& 3160 &  3403  &4\\
38& 584 & 276 &741 &2      &83& 3160 &  3486  &4 \\
39& 721 & 292 &780 &2      &84& 3185 &  3570  &4 \\
40& 928 & 315 &820 &2      &85& 3294 &  3655  &4\\
41& &341 & 861 &2          &86& 3408 &  3741  &4 \\
42& &370 & 903  &2         &87& 3522 &  3828  &4 \\
43& &422 & 946  &2         &88& 3645 &  3916  &4 \\
44& &540 & 990 &2          &89& 3749 &  4005  &4 \\
45& &736 & 1035  &2        &90& 3905 &  4095  &4 \\
$46^\ast$& &1127& 1081  &2       &91&  4038&  4186  &4 \\
47& &1128 & 1128  &2       &92& 4171&  4278  &4\\
48& &1128 & 1176  &2       &93& 4335&  4371  &4 \\
49& &1128  &1225 &2         &$94^\ast$& 4492&  4465 &4\\
50& &1128 & 1275  &4       &$95^\ast$& 4668&  4560  &4\\
51& &1128 & 1326 &4        &$96^\ast$& 4828&  4656  &4\\
\hline
\end{tabular}
}
\end{table}

\end{document}